\documentclass[11pt]{article}
\usepackage{amsmath}
\usepackage{graphicx}
\usepackage{epsfig}
\usepackage{multirow}
\usepackage{color}
\def\disp{\displaystyle}

\def\tto{\;{\lower 1pt \hbox{$\rightarrow$}}\kern -10pt
\hbox{\raise 2pt \hbox{$\rightarrow$}}\;}
\def\Hat{\widehat}

\def\Bar{\overline}
\def\ra{\rangle}
\def\la{\langle}

\def\B{I\!\!B}
\def\h{\hfill\Box}
\def\R{I\!\!R}

\def\ox{\bar{x}}

\def\dom{\mbox{\rm dom}\,}

\def\h{\hfill\triangle}

\def\ph{\varphi}
\def\emp{\emptyset}

\def\oR{\Bar{\R}}
\def\lm{\lambda}
\def\gg{\gamma}
\def\dd{\delta}
\def\al{\alpha}

\newcounter{lk}

\begin{document}
\begin{center}
\vspace*{0.3in} {\bf SOLVING A GENERALIZED HERON PROBLEM\\ BY MEANS OF CONVEX ANALYSIS}\\[2ex]
BORIS S. MORDUKHOVICH\footnote{Department of Mathematics, Wayne
State University, Detroit, MI 48202, USA (email:
boris@math.wayne.edu). Research of this author was partially
supported by the US National Science Foundation under grants
DMS-0603846 and DMS-1007132 and by the Australian Research Council
under grant DP-12092508.}, NGUYEN MAU NAM\footnote{Department of
Mathematics, The University of Texas--Pan American, Edinburg, TX
78539--2999, USA (email: nguyenmn@utpa.edu).} and JUAN
SALINAS\footnote{Department of Mathematics, The University of
Texas--Pan American, Edinburg, TX 78539--2999, USA (email:
jsalinasn@broncs.utpa.edu).}\\[2ex]
\end{center}
\small{\bf Abstract} The classical Heron problem states: \emph{on a
given straight line in the plane, find a point $C$ such that the sum
of the distances from $C$ to the given points $A$ and $B$ is
minimal}. This problem can be solved using standard geometry or
differential calculus. In the light of modern convex analysis, we
are able to investigate more general versions of this problem. In
this paper we propose and solve the following problem: on a given
nonempty closed convex subset of $\R^s$, find a point such that the
sum of the distances from that point to $n$ given nonempty closed
convex subsets of $\R^s$ is minimal.

\newtheorem{Theorem}{Theorem}[section]
\newtheorem{Proposition}[Theorem]{Proposition}
\newtheorem{Remark}[Theorem]{Remark}
\newtheorem{Lemma}[Theorem]{Lemma}
\newtheorem{Corollary}[Theorem]{Corollary}
\newtheorem{Definition}[Theorem]{Definition}
\newtheorem{Example}[Theorem]{Example}
\renewcommand{\theequation}{\thesection.\arabic{equation}}
\normalsize

\section{Problem Formulation}
\setcounter{equation}{0}

Heron from Alexandria (10--75 AD) was ``a Greek geometer and
inventor whose writings preserved for posterity a knowledge of the
mathematics and engineering of Babylonia, ancient Egypt, and the
Greco-Roman world" (from the Encyclopedia Britannica). One of the
geometric problems he proposed in his {\em Catroptica} was as
follows: find a point on a straight line in the plane such that the
sum of the distances from it to two given points is
minimal.\vspace*{0.05in}

Recall that a subset $\Omega$ of $\R^s$ is called {\em convex} if
$\lambda x+(1-\lambda) y\in\Omega$ whenever $x$ and $y$ are in
$\Omega$ and $0\le\lambda\le 1$. Our idea now is to consider a much
broader situation, where two given points in the classical Heron
problem are replaced by finitely many closed and convex subsets
$\Omega_i,\;i=1,\ldots,n$, and the given line is replaced by a given
closed and convex subset $\Omega$ of $\R^s$. We are looking for a
point on the set $\Omega$ such that the sum of the distances from
that point to $\Omega_i,\;i=1,\ldots,n,$ is minimal.\vspace*{0.05in}

The {\em distance} from a point $x$ to a nonempty set $\Omega$ is
understood in the conventional way
\begin{equation}\label{df}
d(x;\Omega)=\inf\big\{||x-y||\;\big|\;y\in\Omega\big\},
\end{equation}
where $||\cdot||$ is the Euclidean norm in $\R^s$. The new {\em
generalized Heron problem} is formulated as follows:
\begin{equation}\label{distance function}
\mbox{minimize }\;D(x):=\sum_{i=1}^n d(x;\Omega_i)\;\mbox{ subject
to } \;x\in \Omega,
\end{equation}
where all the sets $\Omega$ and $\Omega_i$, $i=1,\ldots,n$, are
nonempty, closed, and convex; these are our {\em standing
assumptions} in this paper. Thus \eqref{distance function} is a
constrained convex optimization problem, and hence it is natural to
use techniques of convex analysis and optimization to solve it.

\section{Elements of Convex Analysis}
\setcounter{equation}{0}

In this section we review some basic concepts of convex analysis
used in what follows. This material and much more can be found,
e.g., in the books \cite{bv,HU,r}.

Let $f\colon\R^s\to\oR:=(-\infty,\infty]$ be an extended-real-valued
function, which may be infinite at some points, and let
\begin{equation*}
\mbox{dom }f:=\big\{x\in\R^s\;\big|\;f(x)<\infty\big\}
\end{equation*}
be its {\em effective domain}. The {\em epigraph} of $f$ is a subset
of $\R^s\times\R$ defined by
\begin{equation*}
\mbox{epi }f:=\big\{(x,\alpha)\in\R^{s+1}\;\big|\;x\in\mbox{dom }f
\;\mbox{ and }\;\alpha\ge f(x)\big\}.
\end{equation*}
The function $f$ is {\em closed} if its epigraph is closed, and it
is {\em convex} is its epigraph is a convex subset of $\R^{s+1}$. It
is easy to check that $f$ is convex if and only if
\begin{equation*}
f\big(\lm x+(1-\lm)y\big)\le\lm f(x)+(1-\lm)f(y)\;\mbox{ for all
}\;x, y\in \mbox{dom }f\;\mbox{ and }\;\lm\in[0,1].
\end{equation*}
Furthermore, a nonempty closed subset $\Omega$ of $\R^s$ is convex
if and only if the corresponding distance function
$f(x)=d(x;\Omega)$ is a convex function. Note that the distance
function $f(x)=d(x;\Omega)$ is Lipschitz continuous on $\R^s$ with
modulus one, i.e.,
\begin{equation*}
|f(x)-f(y)|\le||x-y||\;\mbox{ for all }\;x,y\in\R^s.
\end{equation*}
A typical example of an extended-real-valued function is the set
{\em indicator function}
\begin{equation}\label{indicator}
\delta(x;\Omega):=\begin{cases}
0 &\text{if }\;x\in\Omega, \\
\infty & \text{otherwise.}
\end{cases}
\end{equation}
It follows immediately from the definitions that the set
$\Omega\subset\R^s$ is closed (resp.\ convex) if and only if the
indicator function \eqref{indicator} is closed (resp.\ convex).

An element $v\in\R^s$ is called a {\em subgradient} of a convex
function $f\colon\R^s\to\oR$ at $\ox\in\mbox{dom}f$ if it satisfies
the inequality
\begin{equation}\label{convex subdifferential}
\la v, x-\ox\ra\leq f(x)-f(\ox)\;\mbox{ for all }\;x\in\R^s,
\end{equation}
where $\la\cdot,\cdot\ra$ stands for the usual scalar product in
$\R^s$. The set of all the subgradients $v$ in \eqref{convex
subdifferential} is called the {\em subdifferential} of $f$ at $\ox$
and is denoted by $\partial f(\ox)$. If $f$ is convex and
differentiable at $\ox$, then $\partial f(\ox)=\{\nabla f(\ox)\}$.

A well-recognized technique in optimization is to reduce a
constrained optimization problem to an unconstrained one using the
indicator function of the constraint. Indeed, $\ox\in\Omega$ is a
minimizer of the constrained optimization problem:
\begin{equation}\label{c-opt}
\mbox{minimize }\;f(x)\;\mbox{ subject to }\;x\in\Omega
\end{equation}
if and only if it solves the unconstrained problem
\begin{equation}\label{u-opt}
\mbox{minimize }\;f(x)+\delta(x;\Omega),\;\;x\in\R^s.
\end{equation}
By the definitions, for any convex function $\ph\colon\R^s\to\oR$,
\begin{equation}\label{fermat}
\ox\;\mbox{ is a minimizer of }\;\ph\;\mbox{ if and only if
}\;0\in\partial\ph(\ox),
\end{equation}
which is {\em nonsmooth convex} counterpart of the classical {\em
Fermat stationary rule}. Applying \eqref{fermat} to the constrained
optimization problem \eqref{c-opt} via its unconstrained description
\eqref{u-opt} requires the usage of {\em subdifferential calculus}.
The most fundamental calculus result of convex analysis is the
following Moreau-Rockafellar theorem for representing the
subdifferential of sums.

\begin{Theorem}\label{sum rule} Let $\ph_i\colon \R^s\to\oR$,
$i=1,\ldots,m$, be closed convex functions. Assume that there is a
point $\ox\in\cap_{i=1}^n\dom\ph_i$ at which all but $($except
possibly one$)$ of the functions $\ph_1,\ldots,\ph_m$ are
continuous. Then we have the equality
\begin{equation*}
\partial\Big(\sum_{i=1}^m\ph_i\Big)(\ox)=\sum_{i=1}^m\partial\ph_i(\ox).
\end{equation*}
\end{Theorem}

Given a convex set $\Omega\subset\R^s$ and a point $\ox\in\Omega$,
the corresponding geometric counterpart of \eqref{convex
subdifferential} is the {\em normal cone} to $\Omega$ at $\ox$
defined by
\begin{equation}\label{cnc}
N(\ox;\Omega):=\big\{v\in\R^s\big|\;\la v,x-\ox\ra\le 0\;\mbox{ for
all }\;x\in\Omega\big\}.
\end{equation}
It easily follows from the definitions that
\begin{equation}\label{nc}
\partial\delta(\ox;\Omega)=N(\ox;\Omega) \mbox{ for every }\ox\in\Omega,
\end{equation}
which allows us, in particular, to characterize minimizers of the
constrained problem \eqref{c-opt} in terms of the subdifferential
\eqref{convex subdifferential} of $f$ and the normal cone
\eqref{cnc} to $\Omega$ by applying Theorem~\ref{sum rule} to the
function $\ph(x)=f(x)+\dd(x;\Omega)$ in \eqref{fermat}.

Finally in this section, we present a useful formula for computing
the subdifferential of the distance function \eqref{df} via the {\em
unique Euclidean projection}
\begin{eqnarray}\label{pr}
\Pi(\ox;\Omega):=\big\{x\in\Omega\;\big|\;||x-\ox||=d(x;\Omega)\big\}
\end{eqnarray}
of $\ox\in\R^s$ on the closed and convex set $\Omega\subset\R^s$.

\begin{Proposition}\label{subdis} Let $\Omega\ne\emp$ be a closed and
convex of $\R^s$. Then
\begin{eqnarray*}
\partial \mbox{d}(\ox;\Omega)=\left\{\begin{array}{lr}
\Big\{\dfrac{\ox-\Pi(\ox;\Omega)}{d(\ox;\Omega)}\Big\}&\mbox{ if
}\;\ox\notin\Omega,\\\\
N(\ox;\Omega)\cap\B&\mbox{ if }\;\ox\in\Omega,
\end{array}
\right.
\end{eqnarray*}
where $\B$ is the closed unit ball of $\R^s$.
\end{Proposition}

\section{Optimal Solutions to the Generalized Heron Problem}
\setcounter{equation}{0}

In this section we derive efficient characterizations of optimal
solutions to the generalized Heron problem \eqref{distance
function}, which allow us to completely solve this problem in some
important particular settings.\vspace*{0.05in}

First let us present general conditions that ensure the {\em
existence} of optimal solutions to \eqref{distance function}.

\begin{Proposition}\label{ex} Assume that one of the sets $\Omega$ and $\Omega_i$,
$i=1,\ldots,n$, is bounded. Then the generalized Heron problem
\eqref{distance function} admits an optimal solution.
\end{Proposition}
{\bf Proof.} Consider the optimal value
\begin{equation*}
\gamma:=\inf_{x\in\Omega}D(x)
\end{equation*}
in \eqref{distance function} and take a minimizing sequence
$\{x_k\}\subset\Omega$ with $D(x_k)\to\gg$ as $k\to\infty$. If the
constraint set $\Omega$ is bounded, then by the classical
Bolzano-Weierstrass theorem the sequence $\{x_k\}$ contains a
subsequence converging to some point $\ox$, which belongs to the set
$\Omega$ due to it closedness. Since the function $D(x)$ in
\eqref{distance function} is continuous, we have $D(\ox)=\gg$, and
thus $\ox$ is an optimal solution to \eqref{distance function}.

It remains to consider the case when one of sets $\Omega_i$, say
$\Omega_1$, is bounded. In this case we have for the above sequence
$\{x_k\}$ when $k$ is sufficiently large that
\begin{equation*}
d(x_k;\Omega_1)\le D(x_k)<\gamma+1,
\end{equation*}
and thus there exists $w_k\in\Omega_1$ with $||x_k-w_k||<\gamma+1$
for such indexes $k$. Then
\begin{equation*}
||x_k||<||w_k||+\gamma+1,
\end{equation*}
which shows that the sequence $\{x_k\}$ is bounded. The existence of
optimal solutions follows in this case from the arguments above.
$\h$\vspace*{0.05in}

To characterize in what follows optimal solutions to the generalized
Heron problem \eqref{distance function}, for any nonzero vectors
$u,v\in\R^s$ define the quantity
\begin{equation}\label{cos1}
\cos(v,u):=\dfrac{\la v,u\ra}{||v||\cdot||u||}.
\end{equation}
We say that $\Omega$ has a {\em tangent space} at $\ox$ if there
exists a subspace $L=L(\ox)\ne\{0\}$ such that
\begin{equation}\label{tan}
N(\ox;\Omega)=L^\perp:=\big\{v\in\R^s\big|\;\la v,u\ra=0\;\mbox{
whenever }\;u\in L\big\}.
\end{equation}

The following theorem gives necessary and sufficient conditions for
optimal solutions to \eqref{distance function} in terms of
projections \eqref{pr} on $\Omega_i$ incorporated into quantities
\eqref{cos1}. This theorem and its consequences are also important
in verifying the validity of numerical results in the Section 3.

\begin{Theorem}\label{cos} Consider problem \eqref{distance function}
in which
\begin{equation}\label{emp}
\Omega_i\cap\Omega=\emp\;\mbox{ for all }\;i=1,\ldots,n.
\end{equation}
Given $\ox\in\Omega$, define the vectors
\begin{equation}\label{ai}
a_i(\ox):=\dfrac{\ox-\Pi(\ox;\Omega_i)}{d(\ox;\Omega_i)}\ne 0,\quad
i=1,\ldots,n,
\end{equation}
Then $\ox\in\Omega$ is an optimal solution to the generalized Heron
problem \eqref{distance function} if and only if
\begin{equation}\label{ns}
-\sum_{i=1}^n a_i(\ox)\in N(\ox; \Omega).
\end{equation}
Suppose in addition that the constraint set $\Omega$ has a tangent
space $L$ at $\ox$. Then (\ref{ns}) is equivalent to
\begin{equation}\label{cos2}
\sum_{i=1}^n\cos\big(a_i(\ox),u\big)=0\;\mbox{ whenever }\;u\in
L\setminus\{0\}.
\end{equation}
\end{Theorem}
{\bf Proof.} Fix an optimal solution $\ox$ to problem
\eqref{distance function} and equivalently describe it as an optimal
solution to the following unconstrained optimization problem:
\begin{equation}\label{u}
\mbox{minimize }\;D(x)+\delta(x;\Omega),\quad x\in\R^s.
\end{equation}
Applying the generalized Fermat rule \eqref{fermat} to \eqref{u}, we
characterize $\ox$ by
\begin{equation}\label{f}
0\in\partial\Big(\sum_{i=1}^n
d(\cdot;\Omega_i)+\dd(\cdot;\Omega)\Big)(\ox).
\end{equation}
Since all of the functions $d(\cdot;\Omega_i), \;i=1,\ldots,n$, are
convex and continuous, we employ the subdifferential sum rule of
Theorem~\ref{sum rule} to \eqref{f} and arrive at
\begin{eqnarray}\label{f1}
\begin{array}{ll}
0\in\partial\big(D+\delta(\cdot,\Omega)\big)(\ox)&=\disp\sum_{i=1}^n
\partial d(\ox;\Omega_i)
+N(\ox;\Omega)\\
&=\disp\sum_{i=1}^n a_i(\ox)+ N(\ox;\Omega),
\end{array}
\end{eqnarray}
where the second representation in \eqref{f1} is due to \eqref{nc}
and the subdifferential description of Proposition~\ref{subdis} with
$a_i(\ox)$ defined in \eqref{ai}. It is obvious that (\ref{f1}) and
(\ref{ns}) are equivalent.

Suppose in addition that the constraint set $\Omega$ has a tangent
space $L$ at $\ox$. Then the inclusion \eqref{ns} is equivalent to
\begin{equation*}
0\in \sum_{i=1}^n a_i(\ox)+L^\perp,
\end{equation*}
which in turn can be written in the form
\begin{equation*}
\Big\la\sum_{i=1}^n a_i(\ox),u\Big\ra=0 \;\mbox{ for all } \;u\in L.
\end{equation*}
Taking into account that $||a_i(\ox)||=1$ for all $i=1, \ldots, n$
by \eqref{ai} and assumption \eqref{emp}, the latter equality is
equivalent to
\begin{equation*}
\sum_{i=1}^n \dfrac{\la a_i(\ox), v\ra}{||a_i(\ox)||\cdot||u||}=0
\;\mbox{ for all }\;u\in L\setminus\{0\},
\end{equation*}
which gives \eqref{cos2} due to the notation \eqref{cos1} and thus
completes the proof of the theorem. $\h$ \vspace*{0.05in}

To further specify the characterization of Theorem~\ref{cos}, recall
that a set $A$ of $\R^s$ is an {\em affine subspace} if there is a
vector $a\in A$ and a (linear) subspace $L$ such that $A=a+L$. In
this case we say that $A$ is parallel to $L$. Note that the subspace
$L$ parallel to $A$ is uniquely defined by $L=A-A=\{x-y\;|\;x\in A,
\;y\in A\}$ and that $A=b+L$ for any vector $b\in A$.

\begin{Corollary} Let $\Omega$ be an affine subspace parallel to a subspace $L$,
and let assumption \eqref{emp} of Theorem~{\rm\ref{cos}} be
satisfied. Then $\ox\in\Omega$ is a solution to the generalized
Heron problem \eqref{distance function} if and only if condition
\eqref{cos2} holds.
\end{Corollary}
{\bf Proof.} To apply Theorem~\ref{cos}, it remains to check that
$L$ is a tangent space of $\Omega$ at $\ox$ in the setting of this
corollary. Indeed, we have $\Omega=\ox+L$, since $\Omega$ is an
affine subspace parallel to $L$. Fix any $v\in N(\ox;\Omega)$ and
get by \eqref{cnc} that $\la v, x-\ox\ra\le 0$ whenever $x\in\Omega$
and hence $\la v,u\ra\le 0$ for all $u\in L$. Since $L$ is a
subspace, the latter implies that $\la v,u\ra=0$ for all $u\in L$,
and thus $N(\ox;\Omega)\subset L^\perp$. The opposite inclusion is
trivial, which gives \eqref{tan} and completes the proof of the
corollary. $\h$\vspace*{0.05in}

The underlying characterization \eqref{cos2} can be easily checked
when the subspace $L$ in Theorem~\ref{cos} is given as a span of
fixed generating vectors.

\begin{Corollary}\label{finite} Let $L={\rm{span}}\{u_1,\ldots,u_m\}$ with
$u_j\ne 0$, $i=1,\ldots,m$, in the setting of
Theorem~{\rm\ref{cos}}. Then $\ox\in\Omega$ is an optimal solution
to the generalized Heron problem \eqref{distance function} if and
only if
\begin{equation}\label{cos rep1}
\sum_{i=1}^n \cos\big(a_i(\ox),u_j\big)=0 \;\mbox{ for all }\;j=1,
\ldots,m.
\end{equation}
\end{Corollary}
{\bf Proof.} We show that \eqref{cos2} is equivalent to \eqref{cos
rep1} in the setting under consideration. Since \eqref{cos2}
obviously implies \eqref{cos rep1}, it remains to justify the
opposite implication. Denote
\begin{equation*}
a:=\sum_{i=1}^n a_i(\ox)
\end{equation*}
and observe that \eqref{cos rep1} yields the condition
\begin{equation}\label{rep2}
\la a,u_j\ra=0\;\mbox{ for all }\;j=1,\ldots m,
\end{equation}
since $u_j\ne 0$ for all $j=1,\ldots,m$ and $||a_i||=1$ for all
$i=1,\ldots,n$. Taking now any vector $u\in L\setminus\{0\}$, we
represent it in the form
\begin{equation*}
u=\sum_{j=1}^m\lambda_j u_j\;\mbox{ with some }\;\lm_j\in\R^n
\end{equation*}
and get from \eqref{rep2} the equalities
\begin{equation*}
\la a,u\ra=\sum_{j=1}^n\lambda_j\la a,u_j\ra=0.
\end{equation*}
This justifies (\ref{cos2}) and completes the proof of the
corollary. $\h$\vspace*{0.05in}

Let us further examine in more detail the case of two sets
$\Omega_1$ and $\Omega_2$ in \eqref{distance function} with the
normal cone to the constraint set $\Omega$ being a straight line
generated by a given vector. This is a direct extension of the
classical Heron problem to the setting when two points are replaced
by closed and convex sets, and the constraint line is replaced by a
closed convex set $\Omega$ with the property above. The next theorem
gives a complete and verifiable solution to the new problem.

\begin{Theorem}\label{two set} Let $\Omega_1$ and $\Omega_2$ be subsets of
$\R^s$ as $s\ge 1$ with $\Omega\cap\Omega_i=\emp$ for $i=1,2$ in
\eqref{distance function}. Suppose also that there is a vector $a\ne
0$ such that $N(\ox;\Omega)={\rm{span}}\{a\}$. The following
assertions hold, where $a_i:=a_i(\ox)$ are defined in \eqref{ai}:

{\bf (i)} If $\ox\in\Omega$ is an optimal solution to
\eqref{distance function}, then
\begin{equation}\label{normalvector}
\mbox{either }\;a_1+a_2=0\;\mbox{ or }\;\cos(a_1,a)=\cos\big(a_2,a).
\end{equation}

{\bf (ii)} Conversely, if $s=2$ and
\begin{equation}\label{normalvector1}
\mbox{either }\;a_1+a_2=0\;\mbox{ or }\;\big[a_1\ne a_2\;\mbox{ and
}\;\cos(a_1,a)=\cos(a_2,a)\big],
\end{equation}
then $\ox\in\Omega$ is an optimal solution to the generalized Heron
problem \eqref{distance function}.
\end{Theorem}
{\bf Proof.} It follows from the above (see the proof of
Theorem~\ref{cos}) that $\ox\in\Omega$ is an optimal solution to
\eqref{distance function} if and only if $-a_1-a_2\in
N(\ox;\Omega)$. By the assumed structure of the normal cone to
$\Omega$ the latter is equivalent to the alternative:
\begin{equation}\label{nv}
\mbox{either }\;a_1+a_2=0 \;\mbox{ or }\;a_1+a_2=\lambda a\;\mbox{
with some }\;\lambda\ne 0.
\end{equation}

To justify (i), let us show that the second equality in \eqref{nv}
implies the corresponding one in \eqref{normalvector}. Indeed, we
have $||a_1||=||a_1||=1$, and thus \eqref{nv} implies that
\begin{align*}
\lambda^2 ||a||^2=||a_1+a_2||^2 =||a_1||^2+||a_2||^2+2\la a_1,
a_2\ra = 2 +2\la a_1,a_2\ra.
\end{align*}
The latter yields in turn that
\begin{align*}
\la a_1,\lambda a\ra&=\la\lambda a-a_2,\lambda a\ra\\
&=\lambda^2||a||^2-\lambda\la a_2,a\ra\\
&=2 +2\la a_1, a_2\ra-\lambda \la a_2,a\ra\\
&=2\la a_2,a_2\ra + 2\la a_1,a_2\ra-\lambda \la a_2,a\ra\\
&=2\la a_2+a_1,a_2\ra-\lambda \la a_2,a\ra\\
&=2\la\lambda a,a_2\ra-\lambda \la a_2,a\ra=\la a_2,\lambda a\ra,
\end{align*}
which ensures that $\la a_1,a\ra=\la a_2,a\ra$ as $\lambda\ne 0$.
This gives us the equality $\cos(a_1, a)=\cos(a_2,a)$ due to
$||a_1||=||a_2||=1$ and $a\ne 0$. Hence we arrive at
\eqref{normalvector}.

To justify (ii), we need to prove that the relationships in
(\ref{normalvector1}) imply the fulfillment of
\begin{equation}\label{nv1}
-a_1-a_2\in N(\ox;\Omega)=\mbox{span}\{a\}.
\end{equation}
If $-a_1-a_2=0$, then \eqref{nv1} is obviously satisfied. Consider
the alternative in(\ref{normalvector1}) when $a_1\ne a_2$ and
$\cos(a_1,a)=\cos(a_2,a)$. Since we are in $\R^2$, represent
$a_1=(x_1,y_1)$, $a_2=(x_2,y_2)$, and $a=(x,y)$ with two real
coordinates. Then by \eqref{cos1} the equality
$\cos(a_1,a)=\cos(a_2,a)$ can be written as
\begin{equation}\label{dotprod}
x_1x+y_1y=x_2x+y_2y,\;\mbox{ i.e., }\;(x_1-x_2)x=(y_2-y_1)y.
\end{equation}
Since $a\ne 0$, assume without loss of generality that $y\ne 0$. By
\begin{equation*}
||a_1||^2=||a_2||^2\Longleftrightarrow x_1^2+y_1^2=x_2^2+y_2^2
\end{equation*}
we have the equality $(x_1-x_2)(x_1+x_2)=(y_2-y_1)(y_2+y_1)$, which
implies by \eqref{dotprod} that
\begin{equation}\label{dp}
y(x_1-x_2)(x_1+x_2)=x(x_1-x_2)(y_2+y_1).
\end{equation}
Note that $x_1\ne x_2$, since otherwise we have from \eqref{dotprod}
that $y_1=y_2$, which contradicts the condition $a_1\ne a_2$ in
\eqref{normalvector1}. Dividing  both sides of \eqref{dp} by
$x_1-x_2$, we get
\begin{equation*}
y(x_1+x_2)=x(y_2+y_1),
\end{equation*}
which implies in turn that
\begin{equation*}
y(a_1+a_2)=y(x_1+x_2,y_1+y_2)=\big(x(y_1+y_2),
y(y_1+y_2)\big)=(y_1+y_2)a.
\end{equation*}
In this way we arrive at the representation
\begin{equation*}
a_1+a_2=\dfrac{y_1+y_2}{y}a
\end{equation*}
showing that inclusion \eqref{nv1} is satisfied. This ensures the
optimality of $\ox$ in \eqref{distance function} and thus completes
the proof of the theorem. $\h$\vspace*{0.05in}

Finally in this section, we present two examples illustrating the
application of Theorem~\ref{cos} and Corollary~\ref{finite},
respectively, to solving the corresponding the generalized and
classical Heron problems.
\begin{Example}{\rm Consider problem (\ref{distance function}) where $n=2$, the sets
$\Omega_1$ and $\Omega_2$ are two point $A$ and $B$ in the plane,
and the constraint $\Omega$ is a disk that does not contain $A$ and
$B$. Condition (\ref{ns}) from Theorem~\ref{cos} characterizes a
solution $M\in\Omega$ to this generalized Heron problem as follows.
In the first case the line segment $AB$ intersects the disk; then
the intersection is a optimal solution. In this case the problem may
actually have infinitely many solutions. Otherwise, there is a
unique point $M$ on the circle such that a \emph{normal vector}
$\overrightarrow n$ to $\Omega$ at $M$ is the angle bisector of
angle $AMB$, and that is the only optimal solution to the
generalized Heron problem under consideration; see Figure~1.}
\end{Example}
\begin{figure}[h]
\hfill
\begin{minipage}[t]{.45\textwidth}
\begin{center}
\epsfig{file=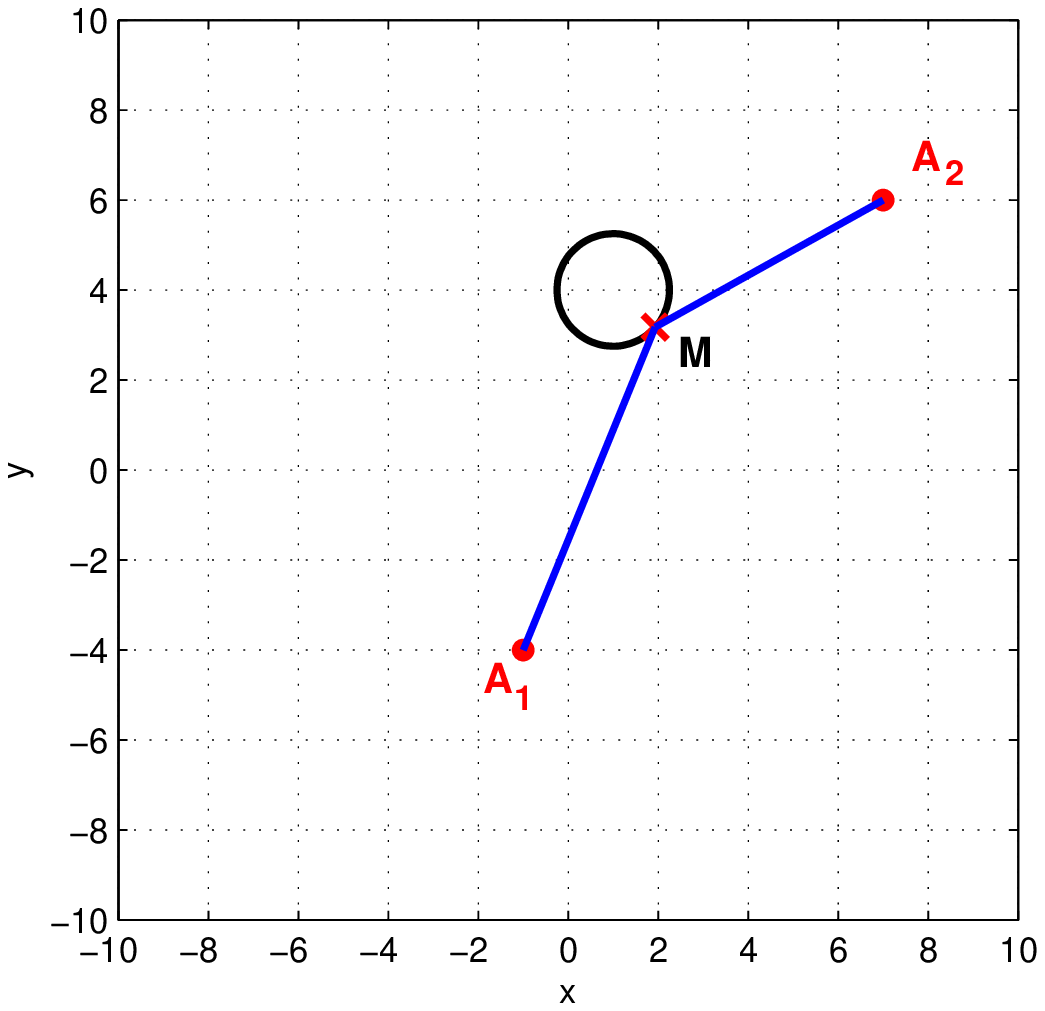, scale=0.5}
\end{center}
\end{minipage}
\hfill
\begin{minipage}[t]{.45\textwidth}
\begin{center}
\epsfig{file=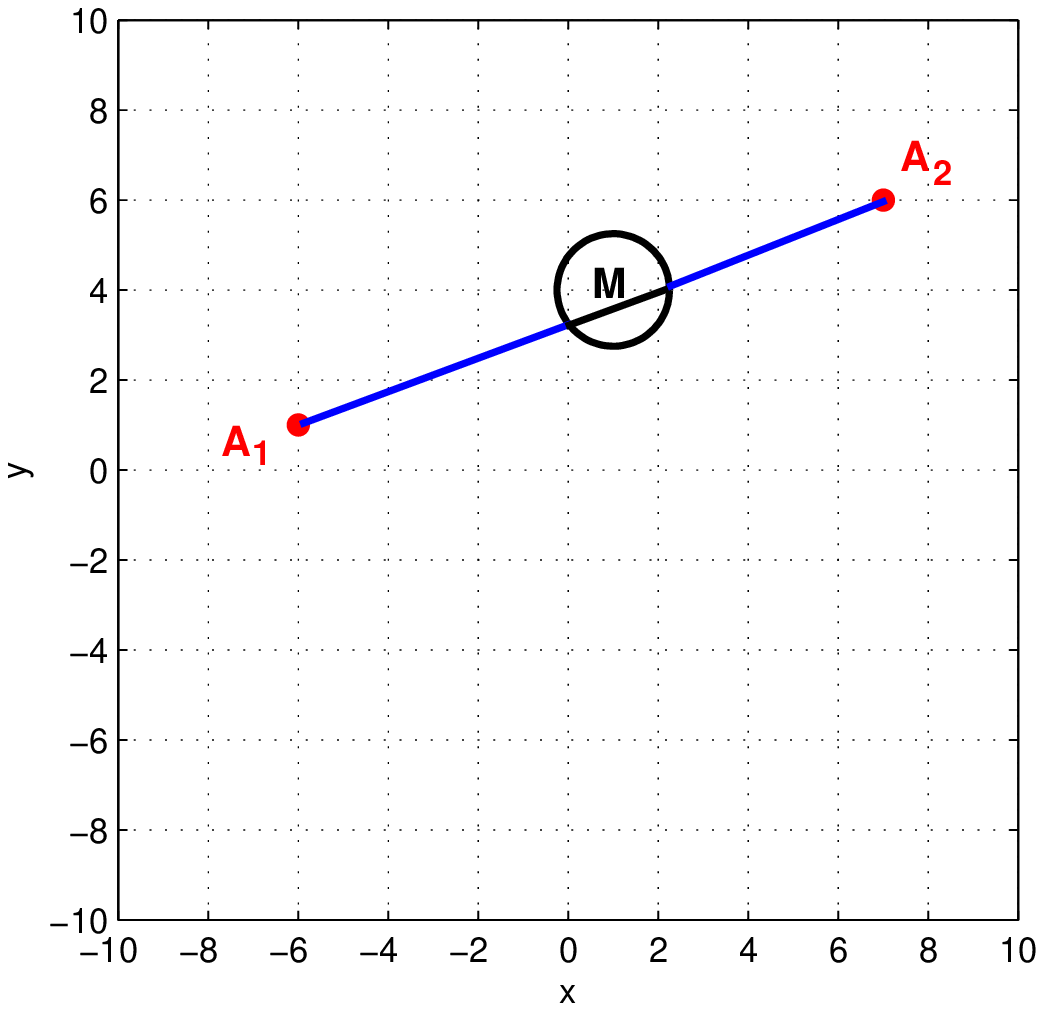, scale=0.5}
\end{center}
\end{minipage}
\hfill\caption{Generalized Heron Problem for Two Points with Disk
Constraint.}
\end{figure}
\begin{Example}{\rm Consider problem (\ref{distance function}),
where $\Omega_i=\{A_i\}$, $i=1,\ldots,n$, are $n$ points in the
plane, and where $\Omega={\cal L}\subset\R^2$ is a straight line
that does not contain these points. Then, by Corollary~\ref{finite}
of Theorem~\ref{cos}, a point $M\in{\cal L}$ is a solution to this
generalized Heron problem if and only if
\begin{equation*}
\cos(\overrightarrow{MA_1},\overrightarrow{a})+\cdots
+\cos(\overrightarrow{MA_n},\overrightarrow{a})=0,
\end{equation*}
where $\overrightarrow{a}$ is a direction vector of ${\cal L}$. Note
that the latter equation completely characterizes the solution of
the classical Heron problem in the plane in both cases when $A_1$
and $A_2$ are on the same side and different sides of ${\cal L}$;
see Figure~2.}
\end{Example}

\begin{figure}[h]
\hfill
\begin{minipage}[t]{.45\textwidth}
\begin{center}
\epsfig{file=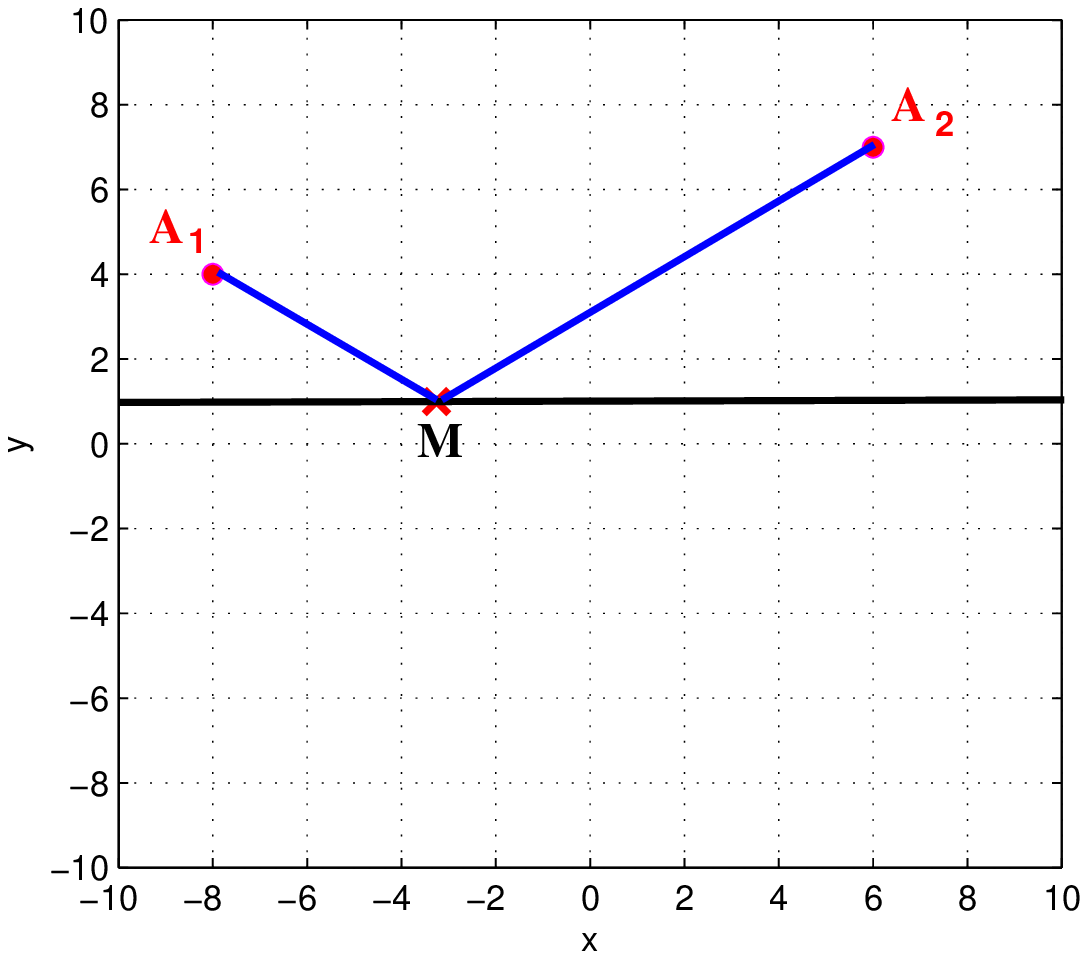, scale=0.5}
\end{center}
\end{minipage}
\hfill
\begin{minipage}[t]{.45\textwidth}
\begin{center}
\epsfig{file=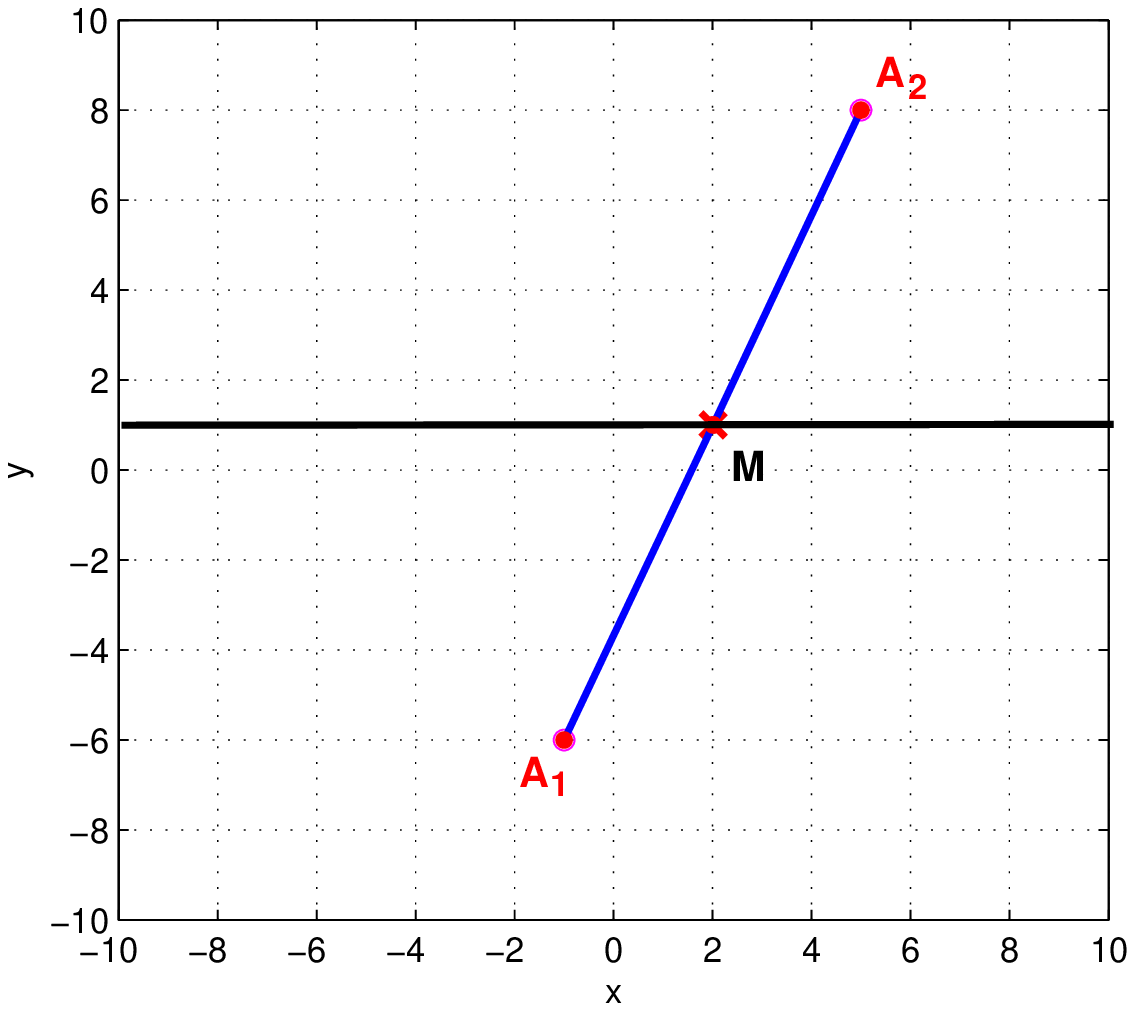, scale=0.5}
\end{center}
\end{minipage}
\hfill \caption{The Classical Heron Problem.}
\end{figure}
\section{Numerical Algorithm and Its Implementation}
\setcounter{equation}{0}

In this section we present and justify an iterative algorithm to
solve the generalized Heron problem \eqref{distance function}
numerically and illustrate its implementations by using MATLAB in
two important settings with disk and cube constraints. Here is the
main algorithm.

\begin{Theorem}\label{subgradient method2} Let $\Omega$ and $\Omega_i$,
$i=1,\ldots,n$, be nonempty closed convex subsets of $\R^s$ such
that at least one of them is bounded. Picking a sequence $\{\al_k\}$
of positive numbers and a starting point $x_1\in \Omega$, consider
the iterative algorithm:
\begin{equation}\label{al}
x_{k+1}=\Pi\Big(x_k-\al_k\sum_{i=1}^n v_{ik};\Omega\Big),\quad
k=1,2,\ldots,
\end{equation}
where the vectors $v_{ik}$ in \eqref{al} are constructed by
\begin{equation}\label{a1}
v_{ik}:=\dfrac{x_k-\omega_{ik}}{d(x_k;\Omega_i)}\;\mbox{ with
}\;\omega_{ik}:=\Pi(x_k;\Omega_i)\;\mbox{ if }\;x_k\notin\Omega_i
\end{equation}
and $v_{ik}:=0$ otherwise. Assume that the given sequence
$\{\al_k\}$ in \eqref{al} satisfies the conditions
\begin{equation}\label{a2}
\sum_{k=1}^\infty\alpha_k=\infty\;\mbox{ and
}\;\sum_{k=1}^\infty\alpha_k^2<\infty.
\end{equation}
Then the iterative sequence $\{x_k\}$ in \eqref{a1} converges to an
optimal solution of the generalized Heron problem \eqref{distance
function} and the value sequence
\begin{equation}\label{Vk}
V_k:=\min\big\{D(x_j)\big|\;j=1,\ldots,k\big\}
\end{equation}
converges to the optimal value $\Hat V$ in this problem.
\end{Theorem}
{\bf Proof.} Observe first of all that algorithm \eqref{al} is well
posed, since the projection to a convex set used in \eqref{a1} is
uniquely defined. Furthermore, all the iterates $\{x_k\}$ in
\eqref{al} are feasible; see the proof of Proposition~\ref{ex}. This
algorithm and its convergence under conditions \eqref{a2} are based
on the subgradient method for convex functions in the so-called
``square summable but not summable case" (see, e.g., \cite{bert}),
the subdifferential sum rule of Theorem~\ref{sum rule}, and the
subdifferential formula for the distance function given in
Proposition~\ref{subdis}. The reader can compare this algorithm and
its justifications with the related developments in \cite{mnft} for
the numerical solution of the (unconstrained) generalized
Fermat-Torricelli problem. $\h$\vspace{0.05in}

Let us illustrate the implementation of the above algorithm and the
corresponding calculations to compute numerically optimal solutions
in the following two characteristic examples.

\begin{Example}\label{disk} {\rm Consider the generalized Heron
problem \eqref{distance function} for pairwise disjoint squares of
{\em right position} in $\R^2$ (i.e., such that the sides of each
square are parallel to the $x$-axis or the $y$-axis) subject to a
given disk constraint. Let $c_i=(a_i,b_i)$ and $r_i$,
$i=1,\ldots,n$, be the centers and the short radii of the squares
under consideration. The vertices of the $i$th square are denoted by
$q_{1i}=(a_i+r_i,b_i+r_i),\;q_{2i}=(a_i-r_i,b_i+r_i),\;q_{3i}=(a_i-r_i,
b_i-r_i),\;q_{4i}=(a_i+r_i,b_i-r_i)$. Let $r$ and $p=(\nu,\eta)$, be
the radius and the center of the constraint. Then the subgradient
algorithm \eqref{al} is written in this case as
\begin{equation*}
x_{k+1}=\Pi\Big(x_k-\al_k\sum_{i=1}^n v_{ik};\Omega\Big),
\end{equation*}
where the projection $P(x,y):=\Pi((x,y);\Omega)$ is calculated by
\begin{equation*}
P(x,y)=(w_x+\nu,w_y+\eta)\;\mbox{ with
}\;w_x=\disp\frac{r(x-\nu)}{\sqrt{(x-\nu)^2+(y-\eta)^2}}\;\mbox{ and
}\; w_y=\disp\frac{r(y-\eta)}{\sqrt{(x-\nu)^2+(y-\eta)^2}}.
\end{equation*}
The quantities $v_{ik}$ in the above algorithm are computed by
{\small\begin{equation*} v_{ik}=\left\{\begin{array}{ll}
0 &\mbox{if }\;|x_{1k}-a_i|\le r_i\;\mbox{ and }\;|x_{2k}-b_i|\le r_i,\\\\
\disp\frac{x_k-q_{1i}}{\|x_k-q_{1i}\|} &\mbox{if }\;x_{1k}-a_i>r_i
\mbox{ and }
\;x_{2k}-b_i> r_i,\\\\
\disp\frac{x_k-q_{2i}}{\|x_k-q_{2i}\|} &\mbox{if
}\;x_{1k}-a_i<-r_i\;
\mbox{ and }\;x_{2k}-b_i> r_i,\\\\
\disp\frac{x_k-q_{3i}}{\|x_k-q_{3i}\|} &\mbox{if
}\;x_{1k}-a_i<-r_i\;
\mbox{ and }\;x_{2k}-b_i< -r_i,\\\\
\disp\frac{x_k-q_{4i}}{\|x_k-q_{4i}\|} &\mbox{if }\;x_{1k}-a_i>r_i\;
\mbox{ and }\;x_{2k}-b_i<-r_i,\\\\
(0,1) &\mbox{if }\;|x_{1k}-a_i|\le r_i \;\mbox{ and }\;x_{2k}-b_i> r_i,\\\\
(0,-1) &\mbox{if }\;|x_{1k}-a_i|\le r_i\;\mbox{ and }\;x_{2k}-b_i< -r_i,\\\\
(1,0) &\mbox{if }\;|x_{1k}-a_i|> r_i\;\mbox{ and }\;|x_{2k}-b_i|\le r_i,\\\\
(-1,0) &\mbox{if }\;|x_{1k}-a_i|< -r_i\;\mbox{ and }\;|x_{2k}-b_i|\le r_i\\\\
\end{array}\right.
\end{equation*}}
for all $i=1,\ldots,n$ and $k=1,2,\ldots$ with the corresponding
quantities $V_k$ defined by \eqref{Vk}.}
\end{Example}

\begin{figure}[h]
\begin{minipage}{2in}
\includegraphics[width=4in]{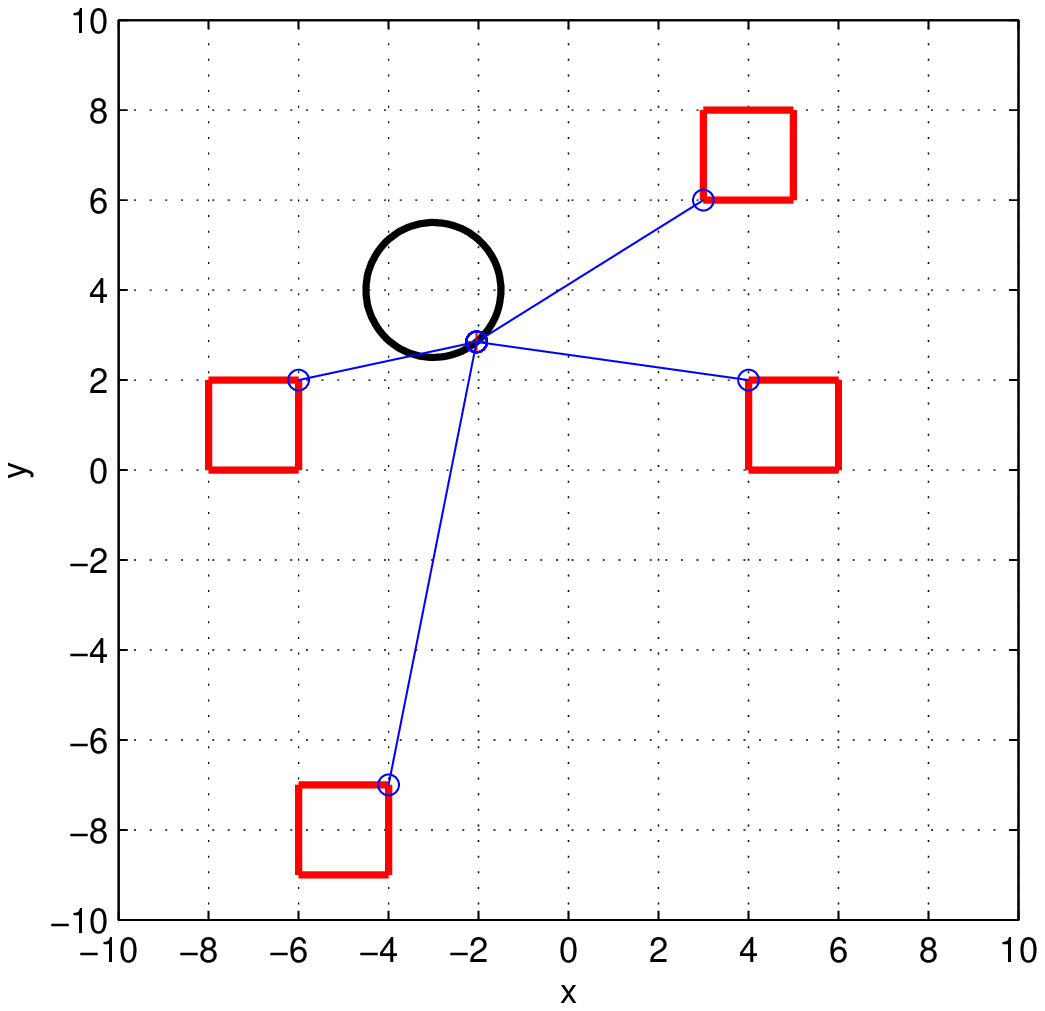}\\
\end{minipage}
~\hfill~
\begin{minipage}[t]{0.52\textwidth}
\begin{tabular}{|c|c|c|}
\hline
\multicolumn{3}{|c|}{MATLAB RESULT} \\
\hline
$k$ & $x_k$ & $V_k$ \\
\hline
1       & (-3,5.5)           & 30.99674 \\
10      & (-1.95277,2.92608) & 26.14035 \\
100     & (-2.02866,2.85698) & 26.13429 \\
1000    & (-2.03861,2.84860) & 26.13419 \\
10,000  & (-2.03992,2.84750) & 26.13419 \\
100,000 & (-2.04010,2.84736) & 26.13419 \\
200,000 & (-2.04011,2.84735) & 26.13419 \\
400,000 & (-2.04012,2.84734) & 26.13419 \\
600,000 & (-2.04012,2.84734) & 26.13419 \\
\hline
\end{tabular}
\end{minipage}
\caption{Generalized Heron Problem for Squares with Disk
Constraint.}
\end{figure}

For the implementation of this algorithm we develop a MATLAB
program. The following calculations are done  and presented below
(see Figure~3 and the corresponding table) for the disk constraint
$\Omega$ with center $(-3,4)$ and radius $1.5$, for the squares
$\Omega_i$ with the same short radius $r=1$ and centers $(-7,1)$,
$(-5,-8)$, (4,7), and (5,1), for the starting point
$x_1=(-3,5.5)\in\Omega$,  and for the sequence of $\al_k=1/k$ in
\eqref{al} satisfying conditions \eqref{a2}. The optimal solution
and optimal value computed up to five significant digits are
$\ox=(-2.04012,2.84734)$ and $\Hat V=26.13419$.\vspace*{0.05in}

The next example concerns the generalized Heron problem for cubes
with ball constraints in $\R^3$.

\begin{Example}\label{cube} {\rm Consider the generalized Heron
problem \eqref{distance function} for pairwise disjoint cubes of
right position in $\R^3$ subject to a ball constraint. In this case
the subgradient algorithm \eqref{al} is
\begin{equation*}
x_{k+1}=\Pi\Big(x_k-\al_k\sum_{i=1}^n v_{ik};\Omega\Big),
\end{equation*}
where the projection $\Pi((x,y,z);\Omega)$ and quantities $v_{ik}$
are computed similarly to Example~\ref{disk}.

\begin{figure}[here]
\begin{minipage}{2in}
   \includegraphics[width=4in]{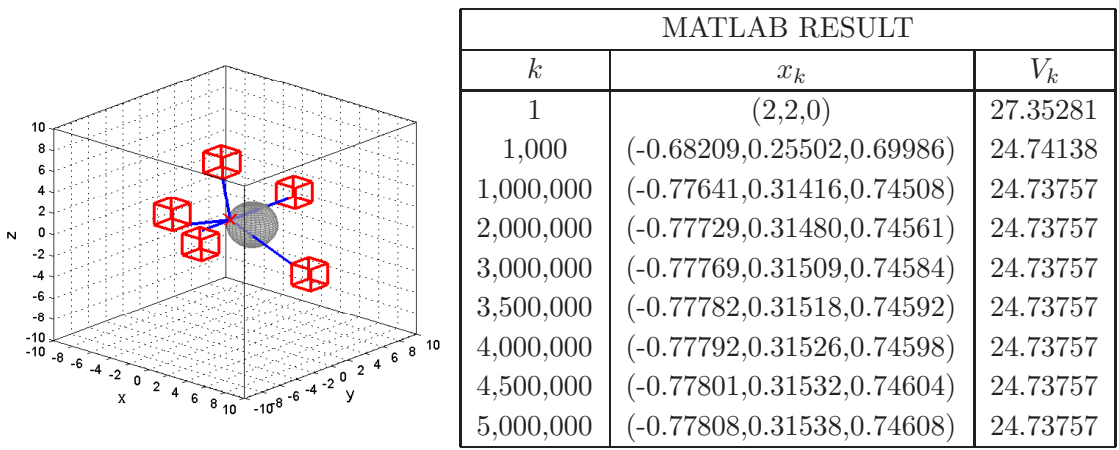}\\
\end{minipage}
~\hfill~
\begin{minipage}[t]{0.56\textwidth}
\begin{tabular}{|c|c|c|}
\hline
\multicolumn{3}{|c|}{MATLAB RESULT} \\
\hline
$k$ & $x_k$ & $V_k$ \\
\hline
1 &         (2,2,0)                    & 27.35281 \\
1,000 &     (-0.68209,0.25502,0.69986) & 24.74138 \\
1,000,000 & (-0.77641,0.31416,0.74508) & 24.73757 \\
2,000,000 & (-0.77729,0.31480,0.74561) & 24.73757 \\
3,000,000 & (-0.77769,0.31509,0.74584) & 24.73757 \\
3,500,000 & (-0.77782,0.31518,0.74592) & 24.73757 \\
4,000,000 & (-0.77792,0.31526,0.74598) & 24.73757 \\
4,500,000 & (-0.77801,0.31532,0.74604) & 24.73757 \\
5,000,000 & (-0.77808,0.31538,0.74608) & 24.73757 \\
\hline
\end{tabular}
\end{minipage}
\caption{Generalized Heron Problem for Cubes with Ball Constraint.}
\end{figure}

For the implementation of this algorithm we develop a MATLAB
program. The Figure~4 and the corresponding figure present the
calculation results for the ball constraint $\Omega$ with center
$(0,2,0)$ and radius $2$, the cubes $\Omega_i$ with centers
$(0,-4,0)$, $(6,2,-3)$, $(-3,-4,2)$, $(-5,4,4)$, and $(-1,8,1)$ with
the same short radius $r=1$, the starting point $x_1=(2,2,0)$, and
the sequence of $\al_k =1/k$ in \eqref{al} satisfying \eqref{a2}.
The optimal solution and optimal value computed up to five
significant digits are $\ox=(-0.77808,0.31538,0.74608)$ and $\Hat
V=24.73756$. }
\end{Example}


\small
\begin{thebibliography}{99}
{\rm \bibitem{bert} Bertsekas, D., Nedic, A., Ozdaglar, A.: Convex
Analysis and Optimization. Athena Scientific, Boston, MA (2003)

\bibitem{bv} Borwein, J.M., Lewis, A.S.: Convex Analysis and
Nonlinear Optimization: Theory and Examples, 2nd edition. Springer,
New York (2006)

\bibitem{HU}Hiriart-Urruty, J.-B., Lemar\'echal, C.: Fundamentals of Convex Analysis. Springer, Berlin (2001)

\bibitem{mnft} Mordukhovich, B.S., Nam, N.M.: Applications of variational analysis to
a generalized Fermat-Torricelli problem. J. Optim. Theory Appl. {\bf
148}, No.\ 3 (2011)

\bibitem{r} Rockafellar, R.T.: Convex Analysis. Princeton University
Press, Princeton, NJ (1970)}
\end{thebibliography}
\end{document}